\documentclass{amsart}
\usepackage{setspace}

\usepackage{amsmath}
\usepackage{amsthm}
\usepackage{amssymb}

\newtheorem{theorem}{Theorem}[section]
\newtheorem{lemma}[theorem]{Lemma}
\newtheorem{proposition}[theorem]{Proposition}
\newtheorem{corollary}[theorem]{Corollary}
\newtheorem{remark}[theorem]{Remark}

\numberwithin{equation}{section}

\newcommand{\Q}{{\mathbb{Q}}}

\newcommand{\OL}{\mathcal{O}}

\begin{document}
\title[On quaternion algebras over the composite of  quadratic number  fields\\
and over some dihedral fields]
{On quaternion algebras over the composite of quadratic  number fields\\
and over some dihedral fields}

\author{Vincenzo Acciaro}
\address{Dipartimento di Economia\\ Universit\`a di Chieti--Pescara\\
Viale Pindaro 42\\  65127 Pescara, Italy}
\email{v.acciaro@unich.it}

\author{Diana Savin}
\address{Faculty of Mathematics and Computer Science, Ovidius University\\
Bd. Mamaia 124, 900527, Constanta, Romania}
\email{savin.diana@univ-ovidius.ro; dianet72@yahoo.com}

\date{}
\begin{abstract}
Let $p$ and $q$ be two positive primes. In this paper we obtain a complete characterization of  quaternion division algebras $H_{K}(p, q)$ over the composite $K$ of $n$ quadratic   number fields. Also, in Section 6, we obtain a characterization of  quaternion division algebras $H_{K}(p, q)$ over some dihedral fields $K.$
\end{abstract}
\maketitle
\section{Introduction} 
\noindent 
Let $F$ be a field with $char(F)\neq2$ and let $a, b \in F \backslash \{0\}$.    {The   generalized quaternion algebra $H_{F}(a, b)$  over the field} $F$ is the algebra having   basis $\{ 1,i,j,k\}$  and     multiplication   table:
\begin{equation*}
\begin{tabular}{c||c|c|c|c|}
$\cdot $ & $1 $ & $i$ & $j$ & $k$ \\ \hline\hline
$1$ & $1$ & $i$ & $j$ & $k$ \\ \hline
$i$ & $i$ & $a$ & $k$ & $a $ \\ \hline
$j$ & $j$ & $-k$ & $b$ & $-b $ \\ \hline
$k$ & $k$ & $-a $ & $b $ & $-ab$ \\ \hline
\end{tabular}%
\end{equation*}%
If  $x=x_{1}  1+ x_{2}i+x_{3}j+x_{4}k \in H_{F}(a, b),$
with $x_{i}\in F$,       the conjugate $\overline{x}$ of $x$ is defined as 
$\overline{x}=x_{1}   1  -x_{2}i -x_{3}j-x_{4}k$, and the norm of $x$ as $\boldsymbol{n}(
x) =x  \overline{x}=x_{1}^{2}-a x_{2}^{2}-b x_{3}^{2}+ ab x_{4}^{2}.$

Quaternion algebras turn out to be  central simple algebras over $F$  (i.e. associative and noncommutative algebras without two sided ideals whose center is precisely $F$)    of dimension $4$ over $F$. 
Recall that the dimension $d$ of  a central simple algebra $A$\ over a  field $F$ is always a perfect square, and its square root $n$ is defined to be   {the degree} of $A$.

The theory of central simple algebras (in particular quaternion algebras and cyclic algebras) has strong connections with algebraic number theory, combinatorics, algebraic geometry, coding theory, computer science and signal theory.

If the equations $ax=b,\ ya=b$   have unique solutions    for all  $a, b \in A, \ a\neq 0$,
then the algebra $A$ is called  {\ a division algebra}.  If $A$ is a
finite-dimensional algebra, then $A$ is a division algebra if and only if $A$
has no zero divisors ($x\neq 0,\  y\neq 0\Rightarrow xy\neq 0$).
In the case of  generalized quaternion algebras  there is a simple criterion that guarantees them to be division algebras:   $H_{F}  (a,b)$ is a division algebra if and only if there is a unique element of zero norm, namely $x = 0$.

Let $L$ be an extension field of $F$, and let $A$ be a central simple algebra over  $F$. 
We recall that $A$ is said to   {split} over $L$,
and $L$ is called
a   {splitting field} for $A$,
if  $A\otimes _{F} L$ is isomorphic with a
matrix algebra over $L$. 

Several results are known about the splitting behavior of quaternion algebras  over specific fields  \cite{gille, lam, voight}.  Explicit conditions which guarantee that a quaternion algebra  splits over the field of rationals numbers -  
or else is a division algebra - were studied in  \cite{alsina}.
In   \cite{savin2016} the second author studied the splitting behavior of some quaternion algebras over some specific quadratic and cyclotomic fields. Moreover, in   \cite{savin2017}  the second author found some sufficient conditions for a quaternion algebra to split over a quadratic field.

In this paper we obtain a complete characterization of  quaternion division algebras $H_{K}(p, q)$ over quadratic and biquadratic number fields $K$, when $p$ and $q$ are two positive primes. 

In this paper, unless otherwise stated, when we say "prime integer" we mean "positive prime integer".

The structure of this paper is the following. In Section 2  we state some results about quaternion algebras and 
quadratic fields  which we will need   later.
In Section 3 we find some necessary and sufficient conditions for a quaternion algebra over a quadratic field to be a division algebra. 
In Section 4 we find some necessary and sufficient conditions for a quaternion algebra over a biquadratic field to be a division algebra. 
In Section 5 we extend the previous results to any composite of $n$ quadratic number fields.
In the last section we compare our approach with the classical ones, from a computational point of view.
\section{ Preliminaries}
In this section we recall some basic facts about quadratic and biquadratic fields, 
as well as some important results concerning quaternion algebras.

Let us recall  first the decomposition behavior of an integral prime ideal in the ring of integers of  quadratic number fields
 \cite[Chapter 13]{ireland}.
\begin{theorem}[Decomposition of primes in quadratic fields]
\label{twodotone}
Let $d\neq 0,1$   {be a square free   integer.}   {Let} $\mathcal{O}_{K}$   {be the ring of integers of the quadratic field} $K=\mathbb{Q}(\sqrt{d})$   {and} $\Delta_{K}$   {be the discriminant of} $K.$   {Let} $p$   {be an odd prime integer. Then, we have}:
\begin{enumerate}
\item $p$   {is ramified in} $\mathcal{O}_{K}$ if and only if
$p  | \Delta_{K}$. In this case 
$p\mathcal{O}_{K}=(p   , \sqrt{d}   )^2$;
\item $p$   splits totally   in $\mathcal{O}_{K}$  if and only if
  $(\frac{\Delta_{K}}{p})=1$. In this case $p\mathcal{O}_{K}=P_1 \cdot  {P_2},$   {where} $P_1$ and $ P_2$   {are distinct prime ideals in} $\mathcal{O}_{K};$
\item $p$   {is inert in} $\mathcal{O}_{K}$ {if and only if}
 $(\frac{\Delta_{K}}{p})=-1$; 
\item 
the prime 2 is ramified in $\mathcal{O}_{K}$  {if and only if}
 $d\equiv 2$ (mod $4$) or $d\equiv 3$ (mod $4$) 
  In the first case    $2\mathcal{O}_{K}=(2,\sqrt{d} )^{2}$, while in the second case
 $2\mathcal{O}_{K}=(2,1+\sqrt{d} )^{2};$
\item   
the prime $2$  splits totally   in $\mathcal{O}_{K}$  {if and only if} 
$d\equiv 1$ (mod 8). In this case     $2\mathcal{O}_{K}=P_1 \cdot P_2$,   {where} $P_1, P_2$
  {are distinct prime ideals in} $\mathcal{O}_{K}$, with $P_1=(2, \frac{1+\sqrt{d}}{2})$;
\item 
the prime $2 $   {is inert  in} $\mathcal{O}_{K}$  {if and only if} 
$d\equiv 5$ (mod $8$). 
\end{enumerate}
\end{theorem}
Now, let $d_{1}$ and $d_{2}$ be two distinct squarefree integers not equal to one.
It is well known that   $K= \mathbb{Q}(\sqrt{d_{1}}, \sqrt{d_{2}}) $ is a Galois extension 
of $\Q$ with Galois group isomorphic to
the Klein $4$-group. 
There are three quadratic subfields of $K$, 
namely $\Q (\sqrt{d_{1}})$,     $\Q (\sqrt{d_{2}} )$ and $\Q  (\sqrt{d_{3} })$, where 
$d_{3}=lcm(  d_{1}, d_{2})/gcd(d_{1},d_{2})$. 
The next result concerning prime ideals which split completely in composita of extensions of the base field
is quite general;  we state it only in the case of our interest, i.e. the biquadratic number fields \cite[p. 46]{marcus}, 
which are the composita of two quadratic number fields.
\begin{theorem}[Splitting of primes in biquadratic fields]
\label{splittingbiquadratic}
Let $d_{1}$ and $d_{2}$ be two distinct squarefree integers not equal to one, and let 
$d_{3}=lcm(  d_{1}, d_{2})/gcd(d_{1},d_{2})$. 
Let $\mathcal{O}_{K}$ denote the ring of integers of the biquadratic field $K$ and
 $\mathcal{O}_{K_{i}}$   the ring of integers of the quadratic subfield $K_{i}=\mathbb{Q}(\sqrt{d_{i}}), \  i=\{1,2,3\}.$
Let $p$ be a prime integer. Then $p$ splits completely in $\OL_K$ if and only if
 $p$ splits completely in each   $\OL_{K_i}, \ (i= \{1,2,3\})$.
 \end{theorem}

Next, let $K$ be a number field  and let $\mathcal{O}_{K}$ be its ring of integers.
If $v$ is a place of $K$, let us denote by $K_v$ the completion of $K$ at $v$.
We recall that a quaternion algebra $H_{K}(a,b)$ is said to ramify
 at a place $v$ of $K$ - or $v$ is said to ramify in $H_{K}(a,b)$ -  if the  quaternion $K_{v}$-algebra 
 $H_{v}=K_{v}\otimes H_{K} (a,b )$
is a division algebra. 
This happens  exactly when the Hilbert symbol
 $(a, b)_{v}$ is equal to $-1$, i.e. when the equation $ax^2+by^2=1$ has no solutions in $K_v$.
We recall that the reduced discriminant $D_{H_{K}(a,b)}$ of the quaternion algebra  $H_{K}(a,b)$
is defined as the product of those prime ideals of the ring of integers $\mathcal{O}_{K}$ of $K$
which  ramify in $H_{K} (a,b )$. 
The following splitting criterion for a quaternion algebras is well  known \cite[Corollary 1.10]{alsina} :
\begin{proposition}
\label{twodotfour}
 Let $K$   be a number field. Then, the quaternion algebra ${H}_{K}(a,b) $  is split if and only if
 its discriminant 
  $D_{H_{K}(a,b)}$ is equal to  the ring of integers $ \mathcal{O}_{K}$ of $K$.
\end{proposition}
If $\mathcal{O}_{K}$ is a principal ideal domain, then we may identify the ideals 
of $\mathcal{O}_{K}$ with their generators, up to units. 
Thus, in a quaternion algebra $H$ over $\mathbb{Q}$, the element $D_{H}$ turns out to be an integer, and $H$  is split if and only if $D_{H}=1$.

The next proposition gives us a geometric interpretation of splitting \cite[Proposition 1.3.2]{gille}:
\begin{proposition}
\label{twodotfive} 
Let  $K$   be a field. Then, the quaternion algebra  ${H}_{K} (a,b ) $ is split if 
and only if the
conic   $C( \alpha ,\beta ) :$ $ax^{2}+by^{2}=z^{2}$ 
 has a rational point in $K$,  i.e.   there are $x_{0},y_{0},z_{0}\in K$   
 such that $ax_{0}^{2}+by_{0}^{2}=z_{0}^{2}$.
\end{proposition}
The next proposition relates the norm group of   extensions of the base 
 field to the splitting behavior of a quaternion  algebra  \cite[ Proposition 1.1.7]{gille}:
\begin{proposition}
\label{twodotsix}  
Let   $F$   be a field. Then, the quaternion algebra 
${H}_{F} (a,b ) $    is split if and only if   $a$ is the norm of an element of $F (\sqrt{b} )$.
 \end{proposition}
For quaternion algebras it is true the following \cite[Proposition 1.1.7]{gille}:
\begin{proposition}
\label{twodoteightone}
Let   $K$   be a field with   char $K \neq 2$   and let  $\alpha, \beta \in K \backslash \{  0 \}$.
 Then the quaternion algebra $H_{K} ( \alpha, \beta )$   is either split or a division algebra.
\end{proposition}
In particular, this tells us that a quaternion algebra $H_\Q (a, b )$  is a division algebra if and only if there is   a prime   
$p$ such that $p | D_{H_\Q  (a, b )}$.
We end up this section with two statements following from the classical 
Albert-Brauer-Hasse-Noether theorem. Proofs of specific formulations 
of this theorem can be found in \cite{linowitz, chinburg}.
\begin{theorem} 
\label{twodotnine}
Let $H_{F}$   be a quaternion algebra over a number field $F$   {and let} $K$   {be a quadratic   extension of} $F.$   {Then there is an embedding of} $K$   {into}
$H_{F}$   {if and only if no prime of} $F$   {which ramifies in} $H_{F}$   splits in  $K.$
\end{theorem}
\begin{proposition} \label{twodotten}
 Let $F$   be a number field and let  $K$   be a quadratic extension of $F$.   Let   $H_{F}$   be a quaternion algebra over  $F$.   Then   $K$   splits $H_{F}$   if and only if there exists an $F$-embedding   $K\hookrightarrow H_{F}$.
\end{proposition}
\section{Division quaternion algebras over   quadratic number fields }
In \cite{savin2016} the second author obtained the following result about quaternion algebras over the field $\mathbb{Q}(i)$:
\begin{proposition}
\label{threedotone}
Let $p\equiv 1$ (mod  $4$)    be a prime 
integer and  let $m$  be an integer which is not a quadratic
residue modulo $p$. Then the quaternion algebra $\mathbb{H}_{\mathbb{Q}( i) }( m ,p) $  is a division
algebra. 
\end{proposition}
In    \cite{savin2017} the second author obtained some sufficient conditions for a quaternion algebra 
$H_{\mathbb{Q}( i) }        (p, q)$ to split, where $p$ and $q$ are two distinct   primes:
\begin{proposition}
\label{threedottwo}
Let $d\neq 0, 1$   {be a   squarefree integer} such that $d\not\equiv 1$ ({mod} $8$),   {and let} $p$ and $q$   {be two primes,} with $q\geq 3$ and $p\neq q.$   {Let} $\mathcal{O}_{K}$   {be the ring of integers of the quadratic field} $K=\mathbb{Q}(\sqrt{d})$   {and} let $\Delta_{K}$   {be the discriminant of} $K$.
\begin{enumerate}
\item
    {if} $p\geq 3$  and  both   $(\frac{\Delta_{K}}{p} ) $ and  $(\frac{\Delta_{K}}{q})$ are not equal to 1,
     {then the quaternion algebra} $H_{\mathbb{Q}(\sqrt{d})}(p,q)$   {splits};
\item   {if} $p=2$   and   $(\frac{\Delta_{K}}{q})\neq 1,$   {then the quaternion algebra} $H_{\mathbb{Q}(\sqrt{d})}(2,q)$   {splits}.
\end{enumerate}
\end{proposition}
From the aforementioned results we   deduce easily  a necessary and sufficient condition for a quaternion algebra $H_{\mathbb{Q}(i)}(p,q)$ to be a division algebra:
\begin{proposition}
\label{threedotthree}
  {Let} $p$ and $ q$   {be two distinct odd primes,}  {such that  } $(\frac{q}{p})\neq 1$.  {Then the quaternion algebra} $H_{\mathbb{Q}(i)}(p,q)$   {is a division algebra if and only if} $p\equiv 1$ ({mod} $4$)   {or} $q\equiv 1$ ({mod} $4$).
\end{proposition}
\begin{proof}
To prove the necessity,  note that if $H_{\mathbb{Q}(i)}(p,q)$ is a division algebra, then   Proposition \ref{threedottwo} and Proposition \ref{twodoteightone} tell us that 
  $(\frac{\Delta_{K}}{p})=1$ or $(\frac{\Delta_{K}}{q})= 1$. This is equivalent to $p\equiv 1$ (mod $4$)   {or} $q\equiv 1$ (mod $4$).
 
To prove the sufficiency, we must distinguish amongst two cases:
\begin{itemize}
\item 
$p\equiv 1$ (mod $4$): \\
Since $(\frac{q}{p})\neq 1$,   Proposition 3.1 tells us that $H_{\mathbb{Q}(i)}(p,q)$ is a division algebra.
\item 
$q\equiv 1$ (mod $4$): \\
Since $(\frac{q}{p})=-1$, the quadratic reciprocity law implies  $(\frac{p}{q})=-1$.
  Proposition 3.1 then tells us that   $H_{\mathbb{Q}(i)}(p,q)$ is a division algebra.
\end{itemize}
\end{proof}
We ask ourselves whether we can obtain a necessary and sufficient explicit condition for    
$H_{\mathbb{Q}(\sqrt{d})}(p,q)$ to be a division algebra when $d$ is arbitrary.
From Proposition \ref{threedottwo} we obtain a necessary explicit condition for   $H_{\mathbb{Q}(\sqrt{d})}(p,q)$ to be a division algebra, namely:   {if} $H_{\mathbb{Q}(\sqrt{d})}(p,q)$   {is a division algebra, then} $(\frac{\Delta_{K}}{p})=1$   {or} $(\frac{\Delta_{K}}{q})=1$. However  this condition is not sufficient:
for example, if we let $K=\mathbb{Q}(\sqrt{3}),\ p=7, \ q=47$, 
then $(\frac{\Delta_{K}}{p})\neq 1$ and $(\frac{\Delta_{K}}{q})=1$,
the quaternion algebra $H_{\mathbb{Q}}(7,47)$ is a division algebra, but
the quaternion algebra $H_{\mathbb{Q}(\sqrt{3})}(7,47)$ is not a division algebra.

It is known   \cite[Section 14.1]{voight} that  if a prime integer $p$ divides  $D_{H(a, b)}$ then it must divide  $2ab$,
hence we may restrict our attention to these primes.
In other words, in order to obtain a sufficient condition for a quaternion algebra 
$H_{\mathbb{Q}(\sqrt{d})}(p,q)$   to be a division algebra,
 it is important to study the ramification of the primes   $ 2,p,q $     in the   algebra    $H_\Q(p,q)$.
The following lemma  \cite[Lemma 1.21]{alsina} gives us a hint:
\begin{lemma}
\label{threedotfour}
Let  $p$ and $q$ be two primes, and let $H_{\mathbb{Q}}(p,q)$ be a quaternion algebra of discriminant $D_H$.
\begin{enumerate}
\item
if $p\equiv q\equiv 3$ (mod  $4$)   {and} $(\frac{q}{p})\neq 1$,   {then} $D_{H}=2p$;
 \item
if $q=2$ and  $p\equiv 3$ ({mod} $8$),  {then} $D_{H}=pq=2p$;
\item
if  $p$   {or} $q\equiv1$ ({mod} $4$), with  $p\neq q$  {and} $(\frac{p}{q})=-1$,   then  $D_{H}=pq$.
\end{enumerate}
\end{lemma}
In addition, the following lemma  \cite[Lemma 1.20]{alsina}     tells us
precisely  when a quaternion algebra $H_{\mathbb{Q}}(p,q)$ splits.
\begin{lemma}
\label{threedotfive}
Let  $p$ and $q$ be two prime integers. Then
$H_{\mathbb{Q}}(p,q)$ is a matrix algebra if and only if one of
the following conditions is satisfied:
\begin{enumerate}
\item
$p=q=2$;
 \item
$p=q \equiv1$ ({mod} $4$);
\item
$q=2$   {and} $p\equiv\pm 1$ ({mod} $8$);
\item
$p\neq q,$ $p\neq 2,$ $q\neq 2,$ $(\frac{q}{p})=1$ and either $p$ or $q$ is congruent to $1$ mod $4.$
\end{enumerate}
\end{lemma}
The next theorem  \cite[Theorem 1.22]{alsina}  describes  the discriminant of
$H_{\mathbb{Q}}(p,q)$, where $p$ and $q$ are primes:
\begin{theorem}
\label{teo}
Let $H = ( \frac{a,b}{\Q} )$ be a quaternion algebra. Then
\begin{enumerate}
\item
If $D_H=1$ then $H$ splits;
\item
if $D_H=2p$, $p$ prime and $p \equiv 3 \pmod{ 4}$ then $H \cong ( \frac{p, -1}{\Q} )$;
\item
if $D_H=pq$, $p,q$ primes, $q \equiv 1 \pmod{4}$ and $(\frac{p}{q})=-1$ then $H \cong ( \frac{p,q}{\Q} )$.
\end{enumerate}
If $a$ and $b$ are prime numbers the algebra $H$ satisfies one and only one of the above statements.
\end{theorem}
We recall that a small ramified $\mathbb{Q}$-algebra is a rational
quaternion algebra having the discriminant equal to the product of two distinct prime numbers.
If we  take into account the previous results and Lemma \ref{threedotfour}, we obtain the following necessary and sufficient explicit condition for a small ramified $\mathbb{Q}$-algebra $H_\Q(p,q)$ to be a division algebra over a quadratic field $\mathbb{Q}(\sqrt{d})$:
\begin{proposition}
\label{threedotsix}
  {Let} $p$ and $q$   {be two distinct  odd primes,}  with  $p$   {or} $q \equiv 1$ ({mod} $4$)   {and} $(\frac{p}{q})=-1.$    Let  $K=\mathbb{Q}(\sqrt{d})$   {and let} $\Delta_{K}$   {be the discriminant of} $K.$
  {Then the quaternion algebra} $H_{\mathbb{Q}(\sqrt{d})}(p,q)$   {is a division algebra if and only if} $(\frac{\Delta_{K}}{p})= 1$   {or} $(\frac{\Delta_{K}}{q})=1$.
\end{proposition}
\begin{proof}
If $H_{\mathbb{Q}(\sqrt{d})}(p,q)$ is a division algebra then  Proposition \ref{threedottwo} gives   $(\frac{\Delta_{K}}{p})=1$ or $(\frac{\Delta_{K}}{q})= 1.$

Conversely, assume that  either $(\frac{\Delta_{K}}{p})= 1$   {or} $(\frac{\Delta_{K}}{q})=1$.
By hypothesis   $p$ or $q$$\equiv$$1$ (mod $4$) and $(\frac{p}{q})=-1.$ 
According to Lemma \ref{threedotfour} (iii), $D_{H}=pq.$ 
This means that the primes which ramify   in the quaternion algebra $H_{\mathbb{Q}}(p,q)$ are precisely
$p$ and $q.$ 
Since either $(\frac{\Delta_{K}}{p})= 1$ or $(\frac{\Delta_{K}}{q})=1$, by  Theorem \ref{twodotone} it follows 
that either $p$ or $q$ splits in the ring of integers of the quadratic field $K.$ 
Finally, Theorem \ref{twodotnine} and Proposition \ref{twodotten} imply that the quaternion algebra $H_{\mathbb{Q}(\sqrt{d})}(p,q)$ does not split, hence, according to Proposition \ref{twodoteightone},
 $H_{\mathbb{Q}(\sqrt{d})}(p,q)$ is a division algebra.
\end{proof}
When  $q=2$ and $p$ is a prime such that $p \equiv 3$ (mod 8),  then,
  according to Lemma \ref{threedotfour} the discriminant 
$D_{H_{\mathbb{Q}(p,q)} }$ is equal to $2p$,    
so $H_{\mathbb{Q}}  (p,q  )$ is a division algebra. 
The next proposition shows what
 happens when we extend the field of scalars from $\Q$ to
$\Q(  \sqrt{d} )$:
\begin{proposition}
\label{threedotseven}  
 {Let} $p$   {be an odd prime,} with $p \equiv 3$ ({mod} $8$).   Let  $K=\mathbb{Q}(\sqrt{d})$   {and let} $\Delta_{K}$   {be the discriminant of} $K.$
  Then  $H_{\mathbb{Q}(\sqrt{d})}(p,2)$   {is a division algebra if and only if} $(\frac{\Delta_{K}}{p})= 1$   {or} $d$ $\equiv 1$ ({mod} $8$).
\end{proposition}
\begin{proof}
 If $H_{\mathbb{Q}(\sqrt{d})}(p,2)$ is a division algebra then, from
  Proposition \ref{threedottwo}, Proposition \ref{twodoteightone}, 
 Theorem \ref{twodotnine} and Proposition \ref{twodotten}, we conclude that $(\frac{\Delta_{K}}{p})=1$ or $d$ $\equiv 1$ (mod $8$). 
 
Conversely, since $p \equiv 3$ (mod $8$) then, according to Lemma \ref{threedotfour}(ii) we must have $D_{H}=2p.$
 It follows that the primes which ramify in   $H_{\mathbb{Q}}(p,2)$ are precisely
$p$ and $2$. Since either $(\frac{\Delta_{K}}{p})= 1$ or $d \equiv 1$ (mod $8$) then, 
after applying Theorem \ref{twodotone}, we obtain that either $p$ or $2$ splits in the ring or integers of $K$. 
From Theorem \ref{twodotnine}, Proposition  \ref{twodotten} and Proposition \ref{twodoteightone}, 
we conclude that  $H_{\mathbb{Q}(\sqrt{d})}(p,2)$ is a division algebra.
\end{proof}
We study next the case where $p$ and $q$ are primes, both congruent to $3$ modulo $4$. 
If $(\frac{q}{p})\neq 1$, then,
according to Lemma \ref{threedotfour}(i), the discriminant 
$D_{H_{\mathbb{Q}}   (p,q)}$ is equal to $2p$, so $H_{\mathbb{Q}}(p,q)$ is a division algebra. 
The next proposition tells us when
the quaternion algebra $H_{\mathbb{Q}(\sqrt{d})}(p,q)$ is still a division algebra.  
\begin{proposition}
\label{threedoteight}
 {Let} $p$ and $q$   {be two odd prime integers,} with $p\equiv q\equiv 3$ ({mod} $4$)  and  $(\frac{q}{p})\neq 1$.   
 Let  $K=\mathbb{Q}(\sqrt{d})$   {and let} $\Delta_{K}$   {be the discriminant of} $K.$
  {Then the quaternion algebra} $H_{\mathbb{Q}(\sqrt{d})}(p,q)$   {is a division algebra if and only if} $(\frac{\Delta_{K}}{p})= 1$   {or} $d$ $\equiv 1$ ({mod} $8$).
\end{proposition}
\begin{proof}
  If $H_{\mathbb{Q}(\sqrt{d})}(p,q)$ is a division algebra then from Proposition \ref{threedottwo}(i) it follows that either
   $(\frac{\Delta_{K}}{p})=1$ or $(\frac{\Delta_{K}}{q})=1$. But, according to Lemma \ref{threedotfour}(i) we must have
    $D_{H_{\mathbb{Q}}(p,q)}=2p.$ So the integral primes   which 
ramify in $H_{\mathbb{Q}}(p,q)$   and could split in $K$ are precisely $p$ and $2.$ 
Finally, after applying Proposition \ref{twodoteightone}, Theorem \ref{twodotnine}, Proposition \ref{twodotten} and Theorem \ref{twodotone},
we obtain that either $(\frac{\Delta_{K}}{p})=1$  or $d$ $\equiv 1$ (mod $8$).

The proof of the converse is similar to  the proof  of sufficiency of Proposition \ref{threedotseven}.
\end{proof}
Taking into account these results and Proposition \ref{twodoteightone}, we are able to understand when a quaternion algebra $H_{\mathbb{Q}(\sqrt{d})}(p,q)$ splits. It is clear that in of the cases  covered by Lemma \ref{threedotfive},
 a quaternion algebra $H_{\mathbb{Q}(\sqrt{d})}(p,q)$ splits. 
Moreover, using Proposition \ref{threedotsix}, Proposition \ref{threedotseven}, 
Proposition \ref{threedoteight} and Proposition \ref{twodoteightone}, we obtain the
following necessary and sufficient explicit condition
for a small ramified $\mathbb{Q}$-algebra $H_\Q(p,q)$ to be a split algebra over a quadratic field $\mathbb{Q}(\sqrt{d})$:
\begin{corollary}
\label{threedotnine}
{Let} $p$ and $q$   {be two distinct  odd primes,}  with  $p$   {or} $q \equiv 1$ ({mod} $4$)   {and} $(\frac{p}{q})=-1$.   
Let      $\Delta_{K}$   {be the discriminant of} $K=\mathbb{Q}(\sqrt{d})$.
  {Then the quaternion algebra} $H_{\mathbb{Q}(\sqrt{d})}(p,q)$ {splits if and only if} $(\frac{\Delta_{K}}{p})\neq1$   {and} $(\frac{\Delta_{K}}{q})\neq1$.
\end{corollary}
\begin{corollary}
\label{threedotten}
{Let} $p$   {be an odd prime,} with $p \equiv 3$ ({mod} $8$).  
Let      $\Delta_{K}$   {be the discriminant of} $K=\mathbb{Q}(\sqrt{d})$.  Then  $H_{\mathbb{Q}(\sqrt{d})}(p,2)$   {splits if and only if} $(\frac{\Delta_{K}}{p})\neq1$   {and} $d$ $\not\equiv 1$ ({mod} $8$).
\end{corollary}
\begin{corollary}
\label{threedoteleven}
{Let} $p$ and $q$   {be two odd prime integers,} with $p\equiv q\equiv 3$ ({mod} $4$)  and  $(\frac{q}{p})\neq 1$.   
 Let      $\Delta_{K}$   {be the discriminant of} $K=\mathbb{Q}(\sqrt{d})$.
  {Then the quaternion algebra} $H_{\mathbb{Q}(\sqrt{d})}(p,q)$   {splits if and only if} $(\frac{\Delta_{K}}{p})\neq1$   {and} $d$  $\not\equiv 1$ ({mod} $8$).
\end{corollary}
The only case left  out is $q = 2, \  p\equiv 5$ ({mod} $8$). 
We consider first the quaternion algebra 
$H_{\mathbb{Q}}(p,q)$, and we get  the following result:
\begin{lemma}
\label{threedotthirteen}
Let  $p\equiv 5$ ({mod} $8$) be a prime integer.
 Then the discriminant of the quaternion algebra $H_{\mathbb{Q}}(p,2 )$ 
is equal to $2p$, and hence $H_{\mathbb{Q}}(p,2)$ is a division algebra.
\end{lemma}
\begin{proof} 
We give here a simple proof which is independent of the theorems stated above.
We know that if a prime divides the discriminant of $H_{\Q}(a,b)$
then it must divide $2ab$.
Since $p\equiv 5$ (mod $8$), from the properties of the Hilbert symbol   and of the
Legendre symbol we obtain:
\[
(2,p)_{p}=(\frac{2}{p})=(-1)^{\frac{p^{2}-1}{8}}=-1
\]
and 
\[
(2,p)_{2}=(-1)^{\frac{p-1}{2}\cdot \frac{1-1}{2} + \frac{p^{2}-1}{8}}=-1
\]
Hence the primes which ramify in $H_{\mathbb{Q}}(p,2)$ are exactly
$2$ and $p$. Therefore,   the reduced discriminant of  $H_{\mathbb{Q}}(p,2 )$  must be equal to $2p$.
\end{proof}
We turn now our attention   to the quaternion algebra $H_{\mathbb{Q}(\sqrt{d})}(p,q)$, where $q = 2$ and $p\equiv 5$ ({mod} $8$).
\begin{proposition}
\label{threedotfourteen}  
 {Let} $p$   {be an odd prime,} with $p \equiv 5$ ({mod} $8$).   Let  $K=\mathbb{Q}(\sqrt{d})$   {and let} $\Delta_{K}$   {be the discriminant of} $K.$ Then  $H_{\mathbb{Q}(\sqrt{d})}(p,2)$   {is a division algebra if and only if} $(\frac{\Delta_{K}}{p})= 1$   {or} $d$ $\equiv 1$ ({mod} $8$).
\end{proposition}
\begin{proof}
The proof is similar to the proof of Proposition \ref{threedotseven},
after replacing   Lemma \ref{threedotfour} with   Lemma \ref{threedotthirteen}. 
\end{proof}
From Proposition \ref{threedotfourteen} and Proposition \ref{twodoteightone} we obtain:
\begin{corollary}
\label{threedotfifteen}
{Let} $p$   {be an odd prime,} with $p \equiv 5$ ({mod} $8$).  
Let      $\Delta_{K}$   {be the discriminant of} $K=\mathbb{Q}(\sqrt{d})$.  Then  $H_{\mathbb{Q}(\sqrt{d})}(p,2)$   {splits if and only if} $(\frac{\Delta_{K}}{p})\neq1$   {and} $d$ $\not\equiv 1$ ({mod} $8$).
\end{corollary}

\begin{theorem}[Classification  over quadratic fields]
\label{maintheorem} 
{Let} $d$  {be a squarefree integer} not equal to one, and let
 $K=\mathbb{Q}(\sqrt{d})$, 
 with discriminant  $\Delta_{K}$. {Let} $p$ and $q$  be two  positive primes. 
 {Then the quaternion algebra} $H_{K}(p,q)$   {is a division algebra if and only if} one of the following conditions
  holds:
 \begin{enumerate}
 \item
$p$ and $q$  are odd and distinct, and \\
$(\frac{p}{q})=-1$,   and \\
     $p \equiv 1 $ ({mod} $4$)  {or} $q  \equiv 1 $ ({mod} $4$),  and  \\
    $(\frac{\Delta_{K}}{p})=1$   {or}  $(\frac{\Delta_{K}}{q})=1$;
 \item
  $q=2$, and  \\
 $p \equiv 3$ ({mod} $8$)   or      $p \equiv 5$ ({mod} $8$), and \\
 either $(\frac{\Delta_{K}}{p})=1$   {or} 
  $d$$\equiv$$1$ ({mod} $8$);
 \item
 $p$ and $q$ are odd,  with $p\equiv q\equiv 3$ ({mod} $4$), and
 \begin{itemize} 
 \item
  $(\frac{q}{p})\neq 1$, and
  \item
 either $(\frac{\Delta_{K}}{p})=1$   {or} $d$$\equiv$$1$ ({mod} $8$);
 \end{itemize}
 or
 \begin{itemize} 
 \item
  $(\frac{p}{q})\neq 1$, and
  \item
 either $(\frac{\Delta_{K}}{q})=1$   {or} $d$$\equiv$$1$ ({mod} $8$);
 \end{itemize}
 \end{enumerate} 
\end{theorem}
\begin{proof}
The theorem follows easily from the four propositions above, after taking into account
Lemma \ref{threedotfour} and Lemma \ref{threedotthirteen}.
\end{proof}
\section{Division quaternion algebras over   biquadratic number fields}
\label{sectionbiquadratic}
Let us start with a  proposition that will be useful   in this section:
\begin{proposition}
\label{fourdotone}
{Let} $a$  {and} $b$  {be distinct nonzero integers and let} $p$ be an odd prime.
\begin{enumerate}
 \item
 \label{uno}
 If $a$ and $b$ are quadratic residues modulo $p,$ then $lcm(a,b)/gcd(a,b)$ is also a quadratic residue modulo $p;$
\item
  \label{due}
If $a$ and $lcm(a,b)/gcd(a,b)$ are quadratic residues modulo $p,$ then $b$ is a quadratic residue modulo $p.$
\end{enumerate}
\end{proposition}
\begin{proof}
Since $a$ and $b$ are both quadratic residues modulo $p$ it follows that   $(\frac{ab}{p})=1$.
Let $c=\frac{lcm(a,b)}{gcd(a,b)}= \frac{ab}{gcd(a,b)^{2}}$.
We have now $(\frac{gcd(a,b)^{2}}{p})=1$ and $(\frac{ab}{p})=(\frac{gcd(a,b)^{2}}{p}) (\frac{c}{p})=1.$ 
Therefore $(\frac{c}{p})=1,$ so $\frac{  lcm(a,b)} {gcd(a,b)}$ is a quadratic residue modulo $p$. 

The proof  of the second case   is similar to the proof of the first case. 
\end{proof}
Let $K/ L$ be an extension of number fields. If a quaternion algebra $H_{L}(p,q)$ splits, then the quaternion algebra $H_{K}(p,q)$ splits as well. 
If the quaternion algebra $H_{L}(p,q)$ is a division algebra then    the quaternion algebra $H_{K}(p,q)$  could   still be
a division algebra or else split.

We consider now a division quaternion algebra $H_{\mathbb{Q}(\sqrt{d_{1}})}(p,q)$  over a base field $L=\mathbb{Q}(\sqrt{d_{1}})$, and   try to find  some conditions which guarantee that it
is still a  division algebra over the biquadratic field $K=\mathbb{Q}(\sqrt{d_{1}}, \sqrt{d_{2}})$. 
 \begin{proposition}
 \label{fourdottwo}
 {Let} $d_{1}$  {and} $d_{2}$  {be distinct   squarefree integers} not equal to one. {Let} $p$ and $q$  {be distinct  odd prime integers} such that $(\frac{p}{q})=-1$, and $p$  {or} $q$ is congruent to one modulo 4.   
    {Let  } $K=\mathbb{Q}(\sqrt{d_{1}}, \sqrt{d_{2}})$  {and let  } $K_{i}=\mathbb{Q}(\sqrt{d_{i}})$ ($i=1,2$),
 with discriminant  $\Delta_{K_{i}}$.
 {Then the quaternion algebra} $H_K(p,q)$   {is a division algebra if and only if} $(\frac{\Delta_{K_{1}}}{p})=(\frac{\Delta_{K_{2}}}{p})=1$   {or} $(\frac{\Delta_{K_{1}}}{q})=(\frac{\Delta_{K_{2}}}{q})=1.$
\end{proposition}
\begin{proof}
Let us assume that $H_K (p,q)$ is a division algebra. 
Hence the quaternion algebras $H_{\mathbb{Q}}(p,q)$,  $H_{K_1}(p,q),$ $H_{K_2}(p,q)$ must  
 all be division algebras.
 
Let $K_{3}=\mathbb{Q}(\sqrt{d_{3}})$ be the third quadratic subfield of $K$, where
$d_{3}= \frac{ lcm(d_{1},d_{2})} {  gcd(d_{1},d_{2})}$,
with discriminant $\Delta_{K_{3}}$. 
Since $H_K (p,q)$ is a division algebra, it follows that $H_{K_3}(p,q)$ must be   a division algebra as well. 
 Proposition \ref{threedotsix} tells us that one of the following conditions must be verified:
\begin{enumerate}
\item
 $(\frac{\Delta_{K_{1}}}{p})=(\frac{\Delta_{K_{2}}}{p})= (\frac{\Delta_{K_{3}}}{p})=1$;
\item
$(\frac{\Delta_{K_{1}}}{q})=(\frac{\Delta_{K_{2}}}{q})= (\frac{\Delta_{K_{3}}}{q})=1$;
\item
There are $i,j$$\in$$\left\{1,2,3\right\},$ $i\neq j$ such that $(\frac{\Delta_{K_{i}}}{p})=(\frac{\Delta_{K_{j}}}{p})=1$ and there is $l$$\in$$\left\{1,2,3\right\},$ $l\neq i,$ $l\neq j$ such that $(\frac{\Delta_{K_{l}}}{q})=1.$ Then $d_{i}$ and $d_{j}$ are quadratic residues modulo $p.$ 
From Proposition \ref{fourdotone} it follows that $d_{1},d_{2}, d_{3}$
are quadratic residues modulo $p,$ so $(\frac{\Delta_{K_{1}}}{p})=(\frac{\Delta_{K_{2}}}{p})= (\frac{\Delta_{K_{3}}}{p})=1$;
\item
There are $i,j$$\in$$\left\{1,2,3\right\},$ $i\neq j$ such that $(\frac{\Delta_{K_{i}}}{q})=(\frac{\Delta_{K_{j}}}{q})=1$ and there is
$l$$\in$$\left\{1,2,3\right\},$ $l\neq i,$ $l\neq j$ such that $(\frac{\Delta_{K_{l}}}{p})=1.$ 
From Proposition \ref{fourdotone} it follows that $(\frac{\Delta_{K_{1}}}{q})=(\frac{\Delta_{K_{2}}}{q})= (\frac{\Delta_{K_{3}}}{q})=1$.
\end{enumerate}
Conversely, let us assume that  $(\frac{\Delta_{K_{1}}}{p})=(\frac{\Delta_{K_{2}}}{p})=1.$ 
By Proposition \ref{fourdotone} it follows that $(\frac{\Delta_{K_{3}}}{p})= 1$.
By Theorem \ref{twodotone}(ii) it follows that $p$ splits in each  $\mathcal{O}_{K_{i}}$. 
According to 
Theorem \ref{splittingbiquadratic},
$p$ must split in $\mathcal{O}_{K}$ as well. 
Since $(\frac{\Delta_{K_{1}}}{p})=1,$    Proposition \ref{threedotsix} tells us that 
$H_{K_1}(p,q)$ must be a division algebra, and that 
$p$ must ramify in $H_{K_1}(p,q).$ 
By Theorem \ref{twodotnine}, Proposition \ref{twodotten} and Proposition \ref{twodoteightone},
  the quaternion algebra  $H_{K}(p,q)$ must be a division algebra.
  
If $(\frac{\Delta_{K_{1}}}{q})=(\frac{\Delta_{K_{2}}}{q})=(\frac{\Delta_{K_{3}}}{q})= 1,$ 
the same argument shows that   $H_{K}(p,q)$ must be a division algebra.
\end{proof}
\begin{proposition}
 \label{fourdotthree}
 {Let} $d_{1}$  {and} $d_{2}$  {be distinct   squarefree integers} not equal to one. {Let} $p$  {an  odd prime integer,} such that $p$$\equiv$$3$ ({mod} $8$). {Let} $K=\mathbb{Q}(\sqrt{d_{1}}, \sqrt{d_{2}})$,  {and let  } $K_{i}=\mathbb{Q}(\sqrt{d_{i}})$ ($i=1,2.$), 
 with discriminant  $\Delta_{K_{i}}$.
 {Then the quaternion algebra} $H_K (p,2)$   {is a division algebra if and only if} $(\frac{\Delta_{K_{1}}}{p})=(\frac{\Delta_{K_{2}}}{p})=1$ {or} $d_{1},$$d_{2}$$\equiv$$1$ ({mod} $8$).
\end{proposition}
\begin{proof} 
Let us assume that $H_K (p,2)$ is a division algebra. 
Then the quaternion algebras $H_{\mathbb{Q}}(p,2),$ $H_{K_1}(p,2),$ $H_{K_2}(p,2)$ are division algebras as well.

Let $K_{3}=\mathbb{Q}(\sqrt{d_{3}})$ be the third quadratic subfield of $K$, where
$d_{3}= \frac{ lcm(d_{1},d_{2})} {  gcd(d_{1},d_{2})}$,
with discriminant $\Delta_{K_{3}}$. 
Since $H_K (p,q)$ is a division algebra, it follows that $H_{K_3}(p,q)$ must be   a division algebra as well.

Let $D_{H}$ be the discriminant of the quaternion algebra $H_{\mathbb{Q}}(p,q).$
According to the hypothesis and to Lemma \ref{threedotfour}(ii), we must have $D_{H}=2p$.
 Proposition \ref{threedotseven} tells us that one of the following conditions must be verified:
\begin{enumerate}
\item
$(\frac{\Delta_{K_{1}}}{p})=(\frac{\Delta_{K_{2}}}{p})= (\frac{\Delta_{K_{3}}}{p})=1$;
\item
 $d_{1},d_{2}, d_{3}$$\equiv$$1$ ({mod} $8$);
\item
There are $i,j$$\in$$\left\{1,2,3\right\},$ $i\neq j$ such that $(\frac{\Delta_{K_{i}}}{p})=(\frac{\Delta_{K_{j}}}{p})=1$ and there is $l$$\in$$\left\{1,2,3\right\},$ $l\neq i,$ $l\neq j$ such that $d_{l}$$\equiv$$1$ ({mod} $8$). Then, $d_{i}$ and $d_{j}$ are quadratic residues modulo $p.$ 
From Proposition \ref{fourdotone} it follows that $d_{1},d_{2}, d_{3}$
are quadratic residues modulo $p,$ so $(\frac{\Delta_{K_{1}}}{p})=(\frac{\Delta_{K_{2}}}{p})= (\frac{\Delta_{K_{3}}}{p})=1$;
\item
There are $i,j$$\in$$\left\{1,2,3\right\},$ $i\neq j$ such that $d_{i},d_{j}$$\equiv$$1$ ({mod} $8$) and there is
$l$$\in$$\left\{1,2,3\right\},$ $l\neq i,$ $l\neq j$ such that $(\frac{\Delta_{K_{l}}}{p})=1.$ 
Since $d_{3}=\frac{lcm(d_{1},d_{2})}{gcd(d_{1},d_{2})}$ it follows    that $d_{1},d_{2}, d_{3}$$\equiv$$1$ ({mod} $8$).
\end{enumerate}
Let us now prove the converse.
If $(\frac{\Delta_{K_{1}}}{p})=(\frac{\Delta_{K_{2}}}{p})= 1$
the argument is the same used in the proof of sufficiency of Proposition \ref{fourdottwo}.

If $d_{1},d_{2}$$\equiv$$1$ ({mod} $8$), it follows easily that  $d_{3}$$\equiv$$1$ ({mod} $8$).
According to Proposition \ref{threedotseven}, the algebras
$H_{K_1}(p,2),$ $H_{K_2}(p,2)$ and $H_{K_3}(p,2)$ are all division algebras and $2$ ramifies there. 
From  Theorem \ref{twodotone} it follows that $2$ splits completely in $K_1$, $K_2$ and $K_3$,
and from Theorem \ref{splittingbiquadratic}
 we can conclude that $2$ splits completely in ${K}.$ 
  Finally, from Theorem \ref{twodotnine}, Proposition \ref{twodotten} and Proposition \ref{twodoteightone},
   it follows now that     $H_{K}(p,2)$ is a division algebra.
\end{proof}
\begin{proposition}
\label{fourdotfour}  
{Let} $d_{1}$  {and} $d_{2}$  {be distinct   squarefree integers} not equal to one. {Let} $p$ and $q$  {be distinct  odd prime integers,}  with $p\equiv q\equiv 3$ ({mod} $4$)  and  $(\frac{q}{p}) = - 1$. {Let} $K=\mathbb{Q}(\sqrt{d_{1}}, \sqrt{d_{2}})$,  {and let  } $K_{i}=\mathbb{Q}(\sqrt{d_{i}})$ ($i=1,2$), with discriminant  $\Delta_{K_{i}}$.
 {Then the quaternion algebra} $H_{K}(p,q)$   {is a division algebra if and only if} $(\frac{\Delta_{K_{1}}}{p})=(\frac{\Delta_{K_{2}}}{p})=1$   {or} $d_{1},d_{2}$$\equiv$$1$ ({mod} $8$).   
 \end{proposition}
\begin{proof}
The proof is similar to the proof of Proposition \ref{fourdotthree}.
\end{proof}
\begin{proposition}
 \label{fourdotfive}
 {Let} $d_{1}$  {and} $d_{2}$  {be distinct   squarefree integers} not equal to one. {Let} $p$  {an  odd prime integer,} such that $p$$\equiv$$5$ ({mod} $8$). {Let} $K=\mathbb{Q}(\sqrt{d_{1}}, \sqrt{d_{2}})$,  {and let  } $K_{i}=\mathbb{Q}(\sqrt{d_{i}})$ ($i=1,2$), 
 with discriminant  $\Delta_{K_{i}}$.
 {Then the quaternion algebra} $H_K (p,2)$   {is a division algebra if and only if} $(\frac{\Delta_{K_{1}}}{p})=(\frac{\Delta_{K_{2}}}{p})=1$ {or} $d_{1},$$d_{2}$$\equiv$$1$ ({mod} $8$).
\end{proposition}
\begin{proof}
The proof is similar to the proof of Proposition \ref{fourdotthree}.
\end{proof}

After gluing the last three proposition together, we obtain the main theorem of our paper:
\begin{theorem}[Classification  over biquadratic fields]
\label{maintheorem} 
{Let} $d_{1}$  {and} $d_{2}$  {be distinct   squarefree integers} not equal to one. {Let} $K=\mathbb{Q}(\sqrt{d_{1}}, \sqrt{d_{2}})$, {and let} $K_{i}=\mathbb{Q}(\sqrt{d_{i}})$ ($i=1,2$) with discriminant  $\Delta_{K_{i}}$. {Let} $p$ and $q$  be two positive
primes. 
 {Then the quaternion algebra} $H_{K}(p,q)$   {is a division algebra if and only if} one of the following conditions
  holds:
 \begin{enumerate}
 \item
$p$ and $q$  are odd and distinct, and \\
$(\frac{p}{q})=-1$,   and \\
     $p \equiv 1 $ ({mod} $4$)  {or} $q  \equiv 1 $ ({mod} $4$),  and  \\
    $(\frac{\Delta_{K_{1}}}{p})=(\frac{\Delta_{K_{2}}}{p})=1$   {or}  $(\frac{\Delta_{K_{1}}}{q})=(\frac{\Delta_{K_{2}}}{q})=1$;
 \item
  $q=2$, and  \\
 $p \equiv 3$ ({mod} $8$)   or      $p \equiv 5$ ({mod} $8$), and \\
 either $(\frac{\Delta_{K_{1}}}{p})=(\frac{\Delta_{K_{2}}}{p})=1$   {or} 
  $d_{1}, d_{2}$$\equiv$$1$ ({mod} $8$);
 \item
 $p$ and $q$ are odd,  with $p\equiv q\equiv 3$ ({mod} $4$), and
 \begin{itemize} 
 \item
  $(\frac{q}{p})\neq 1$, and
  \item
 either $(\frac{\Delta_{K_{1}}}{p})=(\frac{\Delta_{K_{2}}}{p})=1$   {or} $d_{1}, d_{2}$$\equiv$$1$ ({mod} $8$);
 \end{itemize}
 or
 \begin{itemize} 
 \item
  $(\frac{p}{q})\neq 1$, and
  \item
 either $(\frac{\Delta_{K_{1}}}{q})=(\frac{\Delta_{K_{2}}}{q})=1$   {or} $d_{1}, d_{2}$$\equiv$$1$ ({mod} $8$);
 \end{itemize}
 \end{enumerate} 
\end{theorem}
\begin{proof}
The theorem follows easily from the four propositions above, after taking into account
Lemma \ref{threedotfour} and Lemma \ref{threedotthirteen}.
\end{proof}
Let's point out that in the first case of the classification theorem
we have a manifest  symmetry, i.e. by Legendre's statement of the quadratic reciprocity  law $p$
 is a quadratic residue modulo $q$ if and only if $q$ is a quadratic residue mod $p$, while in the third case
the setting is asymmetrical, i.e.,  again by Legendre's statement of the quadratic reciprocity law, $p$
 is a quadratic residue modulo $q$ if and only if $q$ is not a quadratic residue modulo $p$.\\
\section{Extensions}
In this section we show how to apply the technique of proof shown 
in the previous section can be applied   to classify
quaternion algebras over the composite of $n$ quadratic fields.
Let $d_{1},d_{2},..., d_{n}$ be distinct squarefree integers not equal to one.
It is known that   $K= \mathbb{Q}(\sqrt{d_{1}}, \sqrt{d_{2}},..., \sqrt{d_{n}}) $ is a Galois extension 
of $\Q$ with Galois group isomorphic to the group $\underbrace{\mathbb{Z}_{2}\times \mathbb{Z}_{2}\times...\mathbb{Z}_{2}}_{n\  times}.$\\
Let's start with the smallest case, i.e. with $n=3$. For this purpose we take a
quaternion algebra $H_{\mathbb{Q}(\sqrt{d_{1}}, \sqrt{d_{2}})}(p,q)$  
over a base field $L=\mathbb{Q}(\sqrt{d_{1}}, \sqrt{d_{2}})$, and  we try to find  some conditions which guarantee that it
is still a  division algebra over the field $K=\mathbb{Q}(\sqrt{d_{1}}, \sqrt{d_{2}}, \sqrt{d_{3}})$. 

Let us recall the following well known theorem (\cite{waldschmidt}, p.360):
\begin{theorem}
 \label{fourdotseven}
Suppose that $p$ is a prime of $\mathbb{Q}$ which splits completely in each of
two fields $F_{1}$ and $F_{2}$. Then $p$ splits completely in the composite field $F_{1}F_{2}.$ 
\end{theorem}
As a consequence, if $p$ splits completely in a field $F$, then $p$ also splits completely in
the minimal normal extension of $\mathbb{Q}$ containing $F.$
We obtain the following result:
 \begin{proposition}
 \label{fourdoteight}
 {Let} $d_{1}, d_{2}$  {and} $d_{3}$  {be distinct   squarefree integers} not equal to one.
  {Let} $p$ and $q$  {be distinct  odd prime integers} such 
  that $(\frac{p}{q})=-1$, 
  $p$ does not divide $d_{i}$ and $q$ does not divide $d_{i},$ $(i=1,2,3)$ and 
  $p$  {or} $q$ is congruent to one modulo $4.$   
    {Let  } $K=\mathbb{Q}(\sqrt{d_{1}}, \sqrt{d_{2}}, \sqrt{d_{3}})$  {and let  } $K_{i}=\mathbb{Q}(\sqrt{d_{i}})$ ($i=1,2,3$),
 with discriminant  $\Delta_{K_{i}}$.
 {Then the quaternion algebra} $H_K(p,q)$   {is a division algebra if and only if} $(\frac{\Delta_{K_{1}}}{p})=(\frac{\Delta_{K_{2}}}{p})=(\frac{\Delta_{K_{3}}}{p})=1$   {or} $(\frac{\Delta_{K_{1}}}{q})=(\frac{\Delta_{K_{2}}}{q})=(\frac{\Delta_{K_{3}}}{q})=1.$
\end{proposition}
\begin{proof}
If $H_K(p,q)$ is a division algebra, then, according to Proposition \ref{fourdottwo}:
\begin{itemize}
\item 
$H_{\mathbb{Q}(\sqrt{d_{1}}, \sqrt{d_{2}})}(p,q)$ is a division algebra, and this is equivalent to $(\frac{\Delta_{K_{1}}}{p})=(\frac{\Delta_{K_{2}}}{p})=1$ or $(\frac{\Delta_{K_{1}}}{q})=(\frac{\Delta_{K_{2}}}{q})=1;$
\item
$H_{\mathbb{Q}(\sqrt{d_{1}}, \sqrt{d_{3}})}(p,q)$ is a division algebra, and this is equivalent to $(\frac{\Delta_{K_{1}}}{p})=(\frac{\Delta_{K_{3}}}{p})=1$ or $(\frac{\Delta_{K_{1}}}{q})=(\frac{\Delta_{K_{3}}}{q})=1;$ 
\item
$H_{\mathbb{Q}(\sqrt{d_{2}}, \sqrt{d_{3}})}(p,q)$ is a division algebra, and this is equivalent to $(\frac{\Delta_{K_{2}}}{p})=(\frac{\Delta_{K_{3}}}{p})=1$ or $(\frac{\Delta_{K_{2}}}{q})=(\frac{\Delta_{K_{3}}}{q})=1$ 
\end{itemize}
Therefore,  one of the following conditions must be satisfied:
\begin{enumerate}
\item
 $(\frac{\Delta_{K_{1}}}{p})=(\frac{\Delta_{K_{2}}}{p})= (\frac{\Delta_{K_{3}}}{p})=1$;
\item
$(\frac{\Delta_{K_{1}}}{q})=(\frac{\Delta_{K_{2}}}{q})= (\frac{\Delta_{K_{3}}}{q})=1$;
\item
There are $i,j$$\in$$\left\{1,2,3\right\},$ $i\neq j$ such that $(\frac{\Delta_{K_{i}}}{p})=(\frac{\Delta_{K_{j}}}{p})=1$ and there is $l$$\in$$\left\{1,2,3\right\},$ $l\neq i,$ $l\neq j$ such that $(\frac{\Delta_{K_{i}}}{q})=(\frac{\Delta_{K_{l}}}{q})=1$ and $(\frac{\Delta_{K_{j}}}{q})=(\frac{\Delta_{K_{l}}}{q})=1.$ It results that $(\frac{\Delta_{K_{1}}}{q})=(\frac{\Delta_{K_{2}}}{q})= (\frac{\Delta_{K_{3}}}{q})=1.$
\item
There are $i,j,k$$\in$$\left\{1,2,3\right\},$ $i\neq j\neq k\neq i$ such that $(\frac{\Delta_{K_{i}}}{p})=(\frac{\Delta_{K_{j}}}{p})=1$ and $(\frac{\Delta_{K_{j}}}{p})=(\frac{\Delta_{K_{l}}}{p})=1.$ It results that $(\frac{\Delta_{K_{1}}}{p})=(\frac{\Delta_{K_{2}}}{p})= (\frac{\Delta_{K_{3}}}{p})=1.$
\end{enumerate}
Conversely, let us suppose that  $(\frac{\Delta_{K_{1}}}{p})=(\frac{\Delta_{K_{2}}}{p})=(\frac{\Delta_{K_{3}}}{p})=1.$
Since $(\frac{\Delta_{K_{1}}}{p})=(\frac{\Delta_{K_{2}}}{p})=1,$ from   Theorem \ref{splittingbiquadratic} it follows that $p$ splits in $\mathbb{Q}(\sqrt{d_{1}}, \sqrt{d_{2}}).$ Since $(\frac{\Delta_{K_{3}}}{p})=1,$ from Theorem \ref{twodotone}
 it follows that $p$ splits in $\mathbb{Q}(\sqrt{d_{3}}).$
According to Theorem \ref{fourdotseven},   $p$ must split in $K.$
According to Proposition \ref{fourdottwo}, $H_{\mathbb{Q}(\sqrt{d_{1}}, \sqrt{d_{2}})}(p,q)$ is a division algebra and $p$ ramifies in $H_{\mathbb{Q}(\sqrt{d_{1}}, \sqrt{d_{2}})}(p,q).$ By Theorem \ref{twodotnine}, Proposition \ref{twodotten} and Proposition \ref{twodoteightone},
 the quaternion algebra  $H_{K}(p,q)$ must be a division algebra.
  
If $(\frac{\Delta_{K_{1}}}{q})=(\frac{\Delta_{K_{2}}}{q})=(\frac{\Delta_{K_{3}}}{q})= 1,$ 
the same argument shows that   $H_{K}(p,q)$ is a division algebra.
\end{proof}
We  can generalize now Proposition \ref{fourdottwo} and Proposition \ref{fourdoteight}:
\begin{proposition}
 \label{fourdotnine}
 {Let} $n$ {be a positive integer}, $n\geq2$ {and let} $d_{1}, d_{2}...$ $d_{n}$  
 {be distinct   squarefree integers} not equal to one. 
 {Let} $p$ and $q$  {be distinct  odd prime integers} such that $(\frac{p}{q})=-1$, 
 $p$ does not divide $d_{i}$ and $q$ does not divide $d_{i}$ $(i=1,\ldots,n)$, 
 and $p$  {or} $q$ is congruent to one modulo 4.   
    {Let  } $K=\mathbb{Q}(\sqrt{d_{1}}, \sqrt{d_{2}},..., \sqrt{d_{n}})$  {and let  } 
    $K_{i}=\mathbb{Q}(\sqrt{d_{i}}) \ (i=1,\ldots,n)$,
 with discriminant  $\Delta_{K_{i}}$.
 {Then the quaternion algebra} $H_K(p,q)$   {is a division algebra if and only if} $(\frac{\Delta_{K_{1}}}{p})=(\frac{\Delta_{K_{2}}}{p})=...=(\frac{\Delta_{K_{n}}}{p})=1$   {or} $(\frac{\Delta_{K_{1}}}{q})=(\frac{\Delta_{K_{2}}}{q})=...=(\frac{\Delta_{K_{n}}}{q})=1.$
\end{proposition}
\begin{proof} 
The proof is by mathematical induction over $n$, for $n>2$. 
The inductive step is based on the same argument that we  used to go from 
Proposition \ref{fourdottwo} to Proposition \ref{fourdoteight}.
\end{proof}
\begin{proposition}
 \label{fourdotten}
{Let} $d_{1}, d_{2}$  {and} $d_{3}$  {be distinct   squarefree integers} not equal to one. {Let} $p$  {be an odd prime integer} such that $p$ does not divide $d_{i},$ $i=1,2,3,$ $p$$\equiv$$3$ (mod $8$). {Let} $K=\mathbb{Q}(\sqrt{d_{1}}, \sqrt{d_{2}}, \sqrt{d_{3}})$ {and let  } $K_{i}=\mathbb{Q}(\sqrt{d_{i}})$ ($i=1,2,3$), with discriminant  $\Delta_{K_{i}}$.
 {Then the quaternion algebra} $H_K(p,2)$   {is a division algebra if and only if} $(\frac{\Delta_{K_{1}}}{p})=(\frac{\Delta_{K_{2}}}{p})=(\frac{\Delta_{K_{3}}}{p})=1$ {or} $d_{1}, d_{2}, d_{3}$$\equiv$$1$ (mod $8$).
\end{proposition}
\begin{proof}
"$\Rightarrow$" If $H_K(p,2)$ is a division algebra, it results that $H_{\mathbb{Q}(\sqrt{d_{1}}, \sqrt{d_{2}})}(p,q)$ is a division algebra and this is equivalent with $(\frac{\Delta_{K_{1}}}{p})=(\frac{\Delta_{K_{2}}}{p})=1$ or $d_{1}, d_{2}$$\equiv$$1$ (mod $8$);
$H_{\mathbb{Q}(\sqrt{d_{1}}, \sqrt{d_{3}})}(p,q)$ is a division algebra and this is equivalent with $(\frac{\Delta_{K_{1}}}{p})=(\frac{\Delta_{K_{3}}}{p})=1$ or $d_{1}, d_{3}$$\equiv$$1$ (mod $8$); $H_{\mathbb{Q}(\sqrt{d_{2}}, \sqrt{d_{3}})}(p,q)$ is a division algebra and this is equivalent with $(\frac{\Delta_{K_{2}}}{p})=(\frac{\Delta_{K_{3}}}{p})=1$ or $d_{2}, d_{3}$$\equiv$$1$ (mod $8$) (according to Proposition \ref{fourdotthree}). Considering these we can have one of the following cases:
\begin{enumerate}
\item
 $(\frac{\Delta_{K_{1}}}{p})=(\frac{\Delta_{K_{2}}}{p})= (\frac{\Delta_{K_{3}}}{p})=1$;
\item
$d_{1}, d_{2}, d_{3}$$\equiv$$1$ (mod $8$);
\item
There are $i,j$$\in$$\left\{1,2,3\right\},$ $i\neq j$ such that $(\frac{\Delta_{K_{i}}}{p})=(\frac{\Delta_{K_{j}}}{p})=1$ and there is $l$$\in$$\left\{1,2,3\right\},$ $l\neq i,$ $l\neq j$ such that $d_{i}, d_{l}$$\equiv$$1$ (mod $8$) and $d_{j}, d_{l}$$\equiv$$1$ (mod $8$). It results that $d_{1}, d_{2}, d_{3}$$\equiv$$1$ (mod $8$).
\item
There are $i,j,k$$\in$$\left\{1,2,3\right\},$ $i\neq j\neq k\neq i$ such that $(\frac{\Delta_{K_{i}}}{p})=(\frac{\Delta_{K_{j}}}{p})=1$ and $(\frac{\Delta_{K_{i}}}{p})=(\frac{\Delta_{K_{l}}}{p})=1$ and $d_{j}, d_{l}$$\equiv$$1$ (mod $8$). 
It results that $(\frac{\Delta_{K_{1}}}{p})=(\frac{\Delta_{K_{2}}}{p})= (\frac{\Delta_{K_{3}}}{p})=1.$
\end{enumerate}
"$\Leftarrow$" If $(\frac{\Delta_{K_{1}}}{p})=(\frac{\Delta_{K_{2}}}{p})=(\frac{\Delta_{K_{3}}}{p})=1,$ the argument is the same
used in the proof of suffciency of Proposition \ref{fourdoteight}.\\
If $d_{1}, d_{2}, d_{3}$$\equiv$$1$ (mod $8$), applying to Proposition \ref{fourdotthree}, it results that the quaternion algebra  
$H_{\mathbb{Q}(\sqrt{d_{1}}, \sqrt{d_{2}})}(p,q)$ is a division algebras and $2$ ramifies this algebra.
According to Theorem \ref{twodotone} it follows that $2$ splits completely in $K_1$, $K_2$ and $K_3$.
From Theorem \ref{splittingbiquadratic} it results that $2$ splits completely in $\mathbb{Q}(\sqrt{d_{1}}, \sqrt{d_{2}}).$ Using this and the fact that $2$ splits completely in $K_3,$ applying Proposition \ref{fourdotseven}, we conclude that  $2$ splits completely in $K.$ According to  Theorem \ref{twodotnine}, Proposition \ref{twodotten} and Proposition \ref{twodoteightone}, it follows now that     $H_{K}(p,2)$ is a division algebra.
\end{proof}
\begin{proposition}
 \label{fourdoteleven}
{Let} $d_{1}, d_{2}$  {and} $d_{3}$  {be distinct   squarefree integers} not equal to one. {Let} $p,q$  {be two odd prime integers} such that $p$ does not divide $d_{i},$ $i=1,2,3,$ $p$$\equiv$$q$$\equiv$$3$ (mod $4$), $(\frac{q}{p})\neq 1$. {Let} $K=\mathbb{Q}(\sqrt{d_{1}}, \sqrt{d_{2}}, \sqrt{d_{3}})$ {and let  } $K_{i}=\mathbb{Q}(\sqrt{d_{i}})$ ($i=1,2,3$), with discriminant  $\Delta_{K_{i}}$.
 {Then the quaternion algebra} $H_K(p,q)$   {is a division algebra if and only if} $(\frac{\Delta_{K_{1}}}{p})=(\frac{\Delta_{K_{2}}}{p})=(\frac{\Delta_{K_{3}}}{p})=1$ {or} $d_{1}, d_{2}, d_{3}$$\equiv$$1$ (mod $8$).
\end{proposition}
\begin{proof} The proof is similar to the proof of Proposition \ref{fourdotten}.
\end{proof}
\begin{proposition}
 \label{fourdottwelve}
{Let} $d_{1}, d_{2}$  {and} $d_{3}$  {be distinct   squarefree integers} not equal to one. {Let} $p$  {be an odd prime integer} such that $p$ does not divide $d_{i},$ $i=1,2,3,$ $p$$\equiv$$5$ (mod $8$). {Let} $K=\mathbb{Q}(\sqrt{d_{1}}, \sqrt{d_{2}}, \sqrt{d_{3}})$ {and let  } $K_{i}=\mathbb{Q}(\sqrt{d_{i}})$ ($i=1,2,3$), with discriminant  $\Delta_{K_{i}}$.
 {Then the quaternion algebra} $H_K(p,2)$   {is a division algebra if and only if} $(\frac{\Delta_{K_{1}}}{p})=(\frac{\Delta_{K_{2}}}{p})=(\frac{\Delta_{K_{3}}}{p})=1$ {or} $d_{1}, d_{2}, d_{3}$$\equiv$$1$ (mod $8$).
\end{proposition}
\begin{proof}
The proof is similar to the proof of Proposition \ref{fourdotten}.
\end{proof}
Taking into account the results obtained in Proposition \ref{fourdoteight}, Proposition \ref{fourdotten}, Proposition \ref{fourdoteleven}, Proposition \ref{fourdottwelve}, we obtain the following classification theorem.
\begin{theorem}[Classification  over the composite of three quadratic fields]
\label{maintheorem} 
{Let} $d_{1}, d_{2}, d_{3}$  {be distinct   squarefree integers} not equal to one. {Let} $K=\mathbb{Q}(\sqrt{d_{1}}, \sqrt{d_{2}}, \sqrt{d_{3}})$, {and let} $K_{i}=\mathbb{Q}(\sqrt{d_{i}})$ ($i=1,2,3$), with 
discriminant  $\Delta_{K_{i}}$. {Let} $p$ and $q$  be two positive
primes. 
 {Then the quaternion algebra} $H_{K}(p,q)$   {is a division algebra if and only if} one of the following conditions
  holds:
 \begin{enumerate}
 \item
$p$ and $q$  are odd and distinct, and \\
$(\frac{p}{q})=-1$,   and \\
     $p \equiv 1 $ ({mod} $4$)  {or} $q  \equiv 1 $ ({mod} $4$),  and  \\
    $(\frac{\Delta_{K_{1}}}{p})=(\frac{\Delta_{K_{2}}}{p})=(\frac{\Delta_{K_{3}}}{p})=1$   {or}  $(\frac{\Delta_{K_{1}}}{q})=(\frac{\Delta_{K_{2}}}{q})=(\frac{\Delta_{K_{3}}}{q})=1$;
 \item
  $q=2$, and  \\
 $p \equiv 3$ ({mod} $8$)   or      $p \equiv 5$ ({mod} $8$), and \\
 either $(\frac{\Delta_{K_{1}}}{p})=(\frac{\Delta_{K_{2}}}{p})=(\frac{\Delta_{K_{3}}}{p})=1$   {or} 
  $d_{1}, d_{2}, d_{3}$$\equiv$$1$ ({mod} $8$);
 \item
 $p$ and $q$ are odd,  with $p\equiv q\equiv 3$ ({mod} $4$), and
 \begin{itemize} 
 \item
  $(\frac{q}{p})\neq 1$, and
  \item
 either $(\frac{\Delta_{K_{1}}}{p})=(\frac{\Delta_{K_{2}}}{p})=(\frac{\Delta_{K_{3}}}{p})=1$   {or} $d_{1}, d_{2}, d_{3}$$\equiv$$1$ ({mod} $8$);
 \end{itemize}
 or
 \begin{itemize} 
 \item
  $(\frac{p}{q})\neq 1$, and
  \item
 either $(\frac{\Delta_{K_{1}}}{q})=(\frac{\Delta_{K_{2}}}{q})=(\frac{\Delta_{K_{3}}}{q})=1$   {or} $d_{1}, d_{2}, d_{3}$$\equiv$$1$ ({mod} $8$);
 \end{itemize}
 \end{enumerate} 
\end{theorem}
\begin{proof}
The theorem follows easily from the   propositions above, after taking into account
Lemma \ref{threedotfour} and Lemma \ref{threedotthirteen}.
\end{proof}
By mathematical induction over $n$, for $n\geq 2$, we can generalize Proposition \ref{fourdotten}, Proposition \ref{fourdoteleven}, Proposition \ref{fourdottwelve} and   
obtain the following classification for division quaternion algebras over 
a composite of $n$ quadratic fields.
\begin{theorem}[Classification  over a composite of $n$ quadratic fields]
\label{maintheorem} 
{Let}  $d_{1}, d_{2}...$ $d_{n}$  {be distinct   squarefree integers} not equal to one, with $n \geq2$. 
{Let} $p$ and $q$  {be distinct  odd prime integers} such that $p$ does not divide $d_{i}$ and
$q$ does not divide $d_{i}$ $(i= 1,\ldots,n)$. 
{Let  } $K=\mathbb{Q}(\sqrt{d_{1}}, \sqrt{d_{2}},..., \sqrt{d_{n}})$  {and    } 
$K_{i}=\mathbb{Q}(\sqrt{d_{i}})  \   (i= 1,\ldots,n)$,
with discriminant  $\Delta_{K_{i}}$. {Then the quaternion algebra} $H_{K}(p,q)$   {is a division algebra if and only if} one of the following conditions holds:
 \begin{enumerate}
 \item
$p$ and $q$  are odd and distinct, and \\
$(\frac{p}{q})=-1$,   and \\
     $p \equiv 1 $ ({mod} $4$)  {or} $q  \equiv 1 $ ({mod} $4$),  and  \\
    $(\frac{\Delta_{K_{1}}}{p})=(\frac{\Delta_{K_{2}}}{p})=...=(\frac{\Delta_{K_{n}}}{p})=1$   {or}  $(\frac{\Delta_{K_{1}}}{q})=(\frac{\Delta_{K_{2}}}{q})=...=(\frac{\Delta_{K_{n}}}{q})=1$;
 \item
  $q=2$, and  \\
 $p \equiv 3$ ({mod} $8$)   or      $p \equiv 5$ ({mod} $8$), and \\
 either $(\frac{\Delta_{K_{1}}}{p})=(\frac{\Delta_{K_{2}}}{p})=...=(\frac{\Delta_{K_{n}}}{p})=1$   {or} 
  $d_{1}, d_{2}, ..., d_{n}$$\equiv$$1$ ({mod} $8$);
 \item
 $p$ and $q$ are odd,  with $p\equiv q\equiv 3$ ({mod} $4$), and
 \begin{itemize} 
 \item
  $(\frac{q}{p})\neq 1$, and
  \item
 either $(\frac{\Delta_{K_{1}}}{p})=(\frac{\Delta_{K_{2}}}{p})=...=(\frac{\Delta_{K_{n}}}{p})=1$   {or} $d_{1}, d_{2},... ,d_{n}$$\equiv$$1$ ({mod} $8$);
 \end{itemize}
 or
 \begin{itemize} 
 \item
  $(\frac{p}{q})\neq 1$, and
  \item
 either $(\frac{\Delta_{K_{1}}}{q})=(\frac{\Delta_{K_{2}}}{q})=...=(\frac{\Delta_{K_{n}}}{q})=1$   {or} $d_{1}, d_{2},..., d_{n}$$\equiv$$1$ ({mod} $8$).
 \end{itemize}
 \end{enumerate} 
\end{theorem}

\section{Division quaternion algebras over the Hilbert class field of a quadratic field}

We ask ourselves what happens when consider a quaternion algebra over a field $K,$ which is a Galois extension of $\mathbb{Q}$, with Galois group nonabelian, with the order $2l,$ where $l$ is an odd prime integer. For to study this we use the following remark, which can be found in \cite{lam}, p.77.

\begin{remark}
 \label{sixdotone} Let $K/F$ be a finite extension of fields, with the degree $\left[K:F\right]=$odd and let $a, b$$\in$$\dot{F}$ (where $\dot{F}$ is the multiplicative group of the field $F$). Then, the quaternion algebra $H_{K}(a,b)$ splits if and only if $H_{F}(a,b)$ splits.
\end{remark}
First, we study the case when the Galois group Gal($K/\mathbb{Q}$) is isomorphic to the permutations group $S_{3}$ (i.e. the dihedral group $D_{3}$) and we obtain the following results:
\begin{proposition}
 \label{sixdottwo}
{Let} $\epsilon$  {be a primitive root of order} $3$ of the unity and let $F =\mathbb{Q}\left(\epsilon\right)$ be the $3$ th
cyclotomic field. Let $\alpha$$\in$$K\backslash \{0\}$ be a cubicfree integer not equal to one, let the Kummer field $K=F\left(\sqrt[3]{\alpha}\right)$ and {let} $p$ and $q$  {be distinct odd prime integers,} $(\frac{p}{q})=-1$, and $p$  {or} $q$ is congruent to one modulo $4.$ {Then the quaternion algebra} $H_K(p,q)$   {is a division algebra if and only if} $(\frac{-3}{p})=1$ {or} $(\frac{-3}{q})=1.$
\end{proposition}
\begin{proof}
$F =\mathbb{Q}\left(\epsilon\right)=\mathbb{Q}\left(i\sqrt{3}\right)$ and the degree $\left[K:F\right]=3.$\\
$K/\mathbb{Q}$ is a Galois extension and the Galois group Gal($K/\mathbb{Q}$) is isomorphic to the group $S_{3}.$\\
According to Remark \ref{sixdotone} and Proposition \ref{twodoteightone}, $H_K(p,q)$ is a division algebra if and only if $H_F(p,q)$ is a division algebra. Applying to Proposition \ref{threedotsix}, this happens if and only if $(\frac{-3}{p})=1$ {or} $(\frac{-3}{q})=1.$
\end{proof}
\begin{proposition}
 \label{sixdottthree}
{Let} $\epsilon$  {be a primitive root of order} $3$ of the unity and let $F =\mathbb{Q}\left(\epsilon\right)$ be the $3$ th
cyclotomic field. Let $\alpha$$\in$$K\backslash \{0\}$ be a cubicfree integer not equal to one, let the Kummer field $K=F\left(\sqrt[3]{\alpha}\right)$ and {let} $p$ {an odd prime integer,} $p$$\equiv 3$(mod $8$). {Then the quaternion algebra} $H_K(p,2)$   {is a division algebra if and only if} $(\frac{-3}{p})=1.$
\end{proposition}
\begin{proof} The proof is similar with the proof of Proposition \ref{sixdottwo}, but instead of Proposition \ref{threedotsix} we use Proposition \ref{threedotseven}.
\end{proof}
\begin{proposition}
 \label{sixdotfour}
{Let} $\epsilon$  {be a primitive root of order} $3$ of the unity and let $F =\mathbb{Q}\left(\epsilon\right)$ be the $3$ th
cyclotomic field. Let $\alpha$$\in$$K\backslash \{0\}$ be a cubicfree integer not equal to one, let the Kummer field $K=F\left(\sqrt[3]{\alpha}\right)$ and {let} $p$ and $q$  {be distinct odd prime integers,} $(\frac{q}{p})\neq 1$, and $p$$\equiv$$q$$\equiv3$ (mod $4$). {Then the quaternion algebra} $H_K(p,q)$   {is a division algebra if and only if} $(\frac{-3}{p})=1.$
\end{proposition}
\begin{proof}
The proof is similar with the proof of Proposition \ref{sixdottwo}, but instead of Proposition \ref{threedotsix} we use Proposition \ref{threedoteight}.
\end{proof}

Now, we pay attention to the case when the Galois group Gal($K/\mathbb{Q}$) is isomorphic to a dihedral group $D_{l},$ with $l$ prime,
$l\geq 5,$ if this case exists. Here appears the inverse Galois problem. From class field theory, we know that this case exists, that is a dihedral group $D_{l},$ with $l$ prime can be realized as a Galois group over $\mathbb{Q}$ (\cite{jensen}). In \cite{jensen} (p. 352-353) appears the following theorem:
\begin{theorem}
 \label{sixdotfive}
For any prime $l$ and any quadratic field $F =\mathbb{Q}\left(\sqrt{d}\right)$
there exist infinitely many dihedral fields of degree $2l$ containing $F$ (where a dihedral field of degree $2l$ is a normal extension of degree $2l$ over $\mathbb{Q}$ with dihedral Galois group $D_{l}$).
\end{theorem}
Let  $l$ be an odd prime integer and let $F =\mathbb{Q}\left(\sqrt{d}\right)$ be an imaginary quadratic field with class
number $h_{F} =l.$ Let $H_{F}$ be the Hilbert class field of $F.$ If the quaternion algebra $H_{F}(p,q)$ is a division algebra, we are interested when $H_{H_{F}}(p,q)$ when is still a division algebra.
We obtain the following results:
\begin{proposition}
\label{sixdotsix}
{Let} $d$  {be a squarefree integer not equal to one} and let $F =\mathbb{Q}\left(\sqrt{d}\right)$ be an imaginary quadratic field, with class number $h_{F}=l, $ let $H_{F}$ be the Hilbert class field of $F$ and let $\Delta_{F}$ be the discriminant of $F.$ {Let} $p$ and $q$  {be distinct odd prime integers,} $(\frac{p}{q})=-1$, and $p$  {or} $q$ is congruent to one modulo $4.$ {Then the quaternion algebra} $H_{H_{F}}(p,q)$   {is a division algebra if and only if} $\left(\frac{\Delta_{F}}{p}\right)=1$ {or} $\left(\frac{\Delta_{F}}{q}\right)=1.$
\end{proposition}
\begin{proof} The degree $\left[H_{F}:F\right]=l=$odd.\\
It is known that the Hilbert class field $H_{F}$ over $F$ has degree $l,$ $H_{F}/\mathbb{Q}$ is a Galois extension of fields and the Galois group Gal($H_{F}/\mathbb{Q}$) is isomorphic to the dihedral field $D_{l}$ (see \cite{jensen}, p. 348).\\
According to Remark \ref{sixdotone} and Proposition \ref{twodoteightone}, $H_{H_{F}}(p,q)$ is a division algebra if and only if $H_{F}(p,q)$ is a division algebra. Applying to Proposition \ref{threedotsix}, this happens if and only if $\left(\frac{\Delta_{F}}{p}\right)=1$ {or} $\left(\frac{\Delta_{F}}{q}\right)=1.$
\end{proof}
\begin{proposition}
\label{sixdotseven}
{Let} $d$  {be a squarefree integer not equal to one} and let $F =\mathbb{Q}\left(\sqrt{d}\right)$ be an imaginary quadratic field, with class number $h_{F}=l, $ let $H_{F}$ be the Hilbert class field of $F$ and let $\Delta_{F}$ be the discriminant of $F.$ {Let} $p$ {be an odd prime integer,} $p$$\equiv$$3$ (mod $8$). {Then the quaternion algebra} $H_{H_{F}}(p,q)$   {is a division algebra if and only if} $\left(\frac{\Delta_{F}}{p}\right)=1$ {or} $d\equiv1$ (mod $8$).
\end{proposition}
\begin{proof}
The proof is similar with the proof of Proposition \ref{sixdotsix}, when instead of Proposition \ref{threedotsix} we use Proposition \ref{threedotseven}.
\end{proof}
\begin{proposition}
\label{sixdoteight}
{Let} $d$  {be a squarefree integer not equal to one} and let $F =\mathbb{Q}\left(\sqrt{d}\right)$ be an imaginary quadratic field, with class number $h_{F}=l, $ let $H_{F}$ be the Hilbert class field of $F$ and let $\Delta_{F}$ be the discriminant of $F.$ {Let} $p$ and $q$  {be distinct odd prime integers,} $(\frac{q}{p})\neq 1$, and $p$$\equiv$$q$$\equiv$$3$ (mod $4$). {Then the quaternion algebra} $H_{H_{F}}(p,q)$   {is a division algebra if and only if} $\left(\frac{\Delta_{F}}{p}\right)=1$ {or} $d\equiv1$ (mod $8$).
\end{proposition}
\begin{proof}
The proof is similar with the proof of Proposition \ref{sixdotsix}, when instead of Proposition \ref{threedotsix} we use Proposition \ref{threedoteight}.
\end{proof}
Using Theorem \ref{sixdotfive}, we remark that Proposition \ref{sixdotsix}, Proposition \ref{sixdotseven}, Proposition \ref{sixdoteight} remain valid when instead of the Hilbert class field of an imaginary quadratic field $F$ with class number an odd prime $l$ we consider a dihedral field of degree $2l$ containing $F.$

 \section{Final remarks} 
The task of deciding whether a quaternion algebra over a number field is a division algebra
is computationally a feasible one, thanks to the facilities included in   computational algebra packages like
Magma \cite{kohel}, Pari, Sage, etc; however, different approaches lead to very different execution times.

A naive way to check if a quaternion algebra   
is a division algebra is to show that a norm equation has no solution, 
thanks to Proposition \ref{twodotsix} and Proposition \ref{twodoteightone}. 
The problem of determining whether a norm equation over an extension of number fields has a solution
has been extensively investigated in the past, both over arbitrary   and over specific extensions
of  number fields  \cite{acciaro, fieker, fincke, pohst}. 
Algorithms included in  Magma, 
Pari and Sage allow one to find out whether a norm   equation has or not at least
a solution,
and, sometime, to find a sought solution. However in the general case this is not an easy task.

Let us consider first the  apparently efficient   approach based on Proposition \ref{twodotsix}.
For this purpose we wrote two small functions in SAGE, release 8.1, to test the efficiency of this method,
and we used them to test an increasing number of cases (100 - 1000 - 10000 - 100000) of quaternion algebras.

In order to construct the     cases to check,
we considered $r$
unordered couples   $\{p,q\}$ of positive primes $p$ and $q$, for increasing values 
of $p$ and $q$ starting from $2$,
and
\begin{itemize}
\item
in the case $K=\Q (\sqrt{ d })$,  we took
approximately an  equal number $r$    of   squarefree integers  $\pm d$,
for increasing values 
of $d$ starting from one, omitting the trivial case  which would 
give   $K=\Q$;
\item
in the case $K=\Q (\sqrt{ d_1 }   ,   \sqrt{ d_2 }     )$, we took
approximately an  equal number $r$ of
unordered couples $\{\pm d_1, \pm d_2\}$ of   squarefree integers,
for increasing values 
of $d_1$ and $d_2$ starting from $1$, omitting the trivial cases  which would  
give $K=\Q (\sqrt{ t })$ or $K=\Q$.  
\end{itemize}
 
We ran our tests on a  2.6 Ghz - i7 quad core - mac mini (late 2012),  equipped with  8 GB of RAM.
The running time is shown in the second column of
Tables \ref{tavola1} and    \ref{tavola2}.
The value   "n.a" means that the computation was taking too long 
and hence we forced the termination of the program.

\begin{table}[htbp] 
\begin{tabular}{|l|l|l|l|}
\hline 
\# of algebras &  Norm approach & Discriminant  approach & Our approach \\ 
\hline
100 & 1136 & 727 & 3\\
1000 & 16117 & 7941 & 35\\
10000 & n.a. & 82966 & 328\\
100000 & n.a. & 855620 & 3151\\ 
\hline 
\end{tabular}
\caption{Running time in ms. required to test   $H_{\Q(\sqrt{d })}(p,q)$.}\label{tavola1}
\end{table}

A different approach to check if a quaternion algebra   
is a division algebra, based on Proposition \ref{twodotfour}, is to show that the discriminant ideal of the algebra 
is not equal to the full ring of integers of the base field.
Again we wrote two small functions in SAGE to test the efficiency of this approach.
The running time is shown in the third column of
Tables \ref{tavola1} and     \ref{tavola2}.

\begin{table}[htbp] 
\begin{tabular}{|l|l|l|l|}
\hline 
\# of algebras &  Norm approach & Discriminant  approach & Our approach \\ 
\hline
100 & 2501 & 1044 & 4\\
1000 & 43329 & 9428 & 42\\
10000 & n.a. & 97049 & 384\\
100000 & n.a. & 1000279 & 3818\\ 
\hline 
\end{tabular}
\caption{Running time in ms. required to test  $H_{\Q( \sqrt{d_1 },     \sqrt{d_2 }    )}(p,q)$.}\label{tavola2}
\end{table}

The approach described in this paper does not require one to solve norm equations over relative extensions of number fields,
neither to compute the discriminant of  quaternion algebras defined over an arbitrary number field.
In fact, all we need is to compute a few  Legendre symbols as well as 
the discriminants of quadratic extensions of $\Q$ involved, 
which is a very easy  task:
indeed, for a nonzero square free integer  $d$,  the discriminant of the quadratic field $\Q(\sqrt{d})$ 
is $d$ if $d$  is congruent to 1 modulo 4,    otherwise $4d$.
We wrote two small functions in SAGE to test the efficiency of our approach. 
The running time is shown in the fourth column of
Tables \ref{tavola1} and     \ref{tavola2}.

\section*{Acknowledgements} 
Since part of this work   was written  in the summer   2017, while the second author was visiting 
the University  "Gabriele D'Annunzio" of Chieti-Pescara, she wants to thank the   Department of Economic Studies 
for the hospitality and the support. The second author wants to thank Professor Mohammed Taous as well,
for the fruitful discussions about  biquadratic fields.


\begin{thebibliography}{99}

\bibitem{acciaro}
V. Acciaro, Solvability of norm equations over cyclic number fields of prime degree, Mathematics of computation,
vol. 65, No. 216(1996), p. 1663-1674.

\bibitem{alsina}
M. Alsina and P. Bayer, Quaternion Orders, Quadratic Forms and Shimura Curves, CRM Monograph Series,
22, American Mathematical Society, 2004.

\bibitem{magma}
W. Bosma, J. Cannon, and C. Playoust, The Magma algebra system. I. The user language, J. Symbolic Comput., 24 (1997),
pp. 235-265.

\bibitem{chinburg}
T. Chinburg and E. Friedman, An embedding theorem for quaternion algebras, J. London Math. Soc. (2), 60(1):33--44, 1999.

\bibitem{fieker} C. Fieker, A. Jurk and M. Pohst, On solving relative norm equations in algebraic number fields, Math. Computation, 66 (217):399-410, 1997.

\bibitem{fincke} U. Fincke, M. Pohst, A Procedure for Determining Algebraic Integers of Given Norm, Proceedings
EUROCAL 83, Springer Lecture Notes in Computer Science No. 162, Springer, 1983.

\bibitem{gille}
P. Gille, T. Szamuely,  Central Simple Algebras and Galois
Cohomology, Cambridge University Press, 2006.

\bibitem{ireland}
K. Ireland,  M. Rosen,  A Classical Introduction to Modern Number Theory, Springer Verlag, 1992.

\bibitem{jensen} C.U. Jensen, N. Yui, Polynomials with $D_{p}$ as Galois Group, Journal of  Number Theory  {15} (1982), pp. 347-375.

\bibitem{kohel}
D.R. Kohel, Quaternion algebras, available online at 
http://www.i2m.univ-amu.fr/perso/david.kohel/alg/doc/AlgQuat.pdf.


\bibitem{lam}
T. Y. Lam,  Introduction to Quadratic Forms over Fields,
American Mathematical Society, 2004.

\bibitem{ledet}
A. Ledet, Brauer Type Embedding Problems, AMS, 2005.

\bibitem{linowitz}
B. Linowitz,  Selectivity in quaternion algebras, Journal of  Number Theory  {132} (2012), pp. 1425-1437. 

\bibitem{marcus}
D. Marcus,  Number Fields, Springer-Verlag, New York, 1977.

\bibitem{pohst} M. Pohst, H. Zassenhaus, Algorithmic Algebraic Number Theory, Encyclopaedia of Mathematics
and its Applications, Cambridge University Press, 1989. 

\bibitem{savin2016}
D. Savin,  About division quaternion algebras and division symbol algebras,
Carpathian Journal of Mathematics, 32(2), p. 233-240 (2016).

\bibitem{savin2017}
 D. Savin,  About split quaternion algebras over quadratic fields and symbol algebras of degree  $n$, Bull. Math. Soc. Sci. Math. Roumanie, Tome 60 (108) No. 3, 2017, p. 307-312.

\bibitem{voight}
J. Voight,  The Arithmetic of Quaternion Algebras, Available online at
 {http://www.math.dartmouth.edu/{\textasciitilde}jvoight/quat-book.pdf}.

\bibitem{waldschmidt}
M. Waldschmidt, P. Moussa, J.M. Luck, C. Itzykson (editors), From Number Theory to Physics, Springer, 1992 

\end{thebibliography}
\end{document}